\documentstyle[11pt]{article}
\newcommand{\qed}{\hfill$\Box$}
\newcommand{\suml}{\sum\limits}

\newcommand\weight{\hbox{\sl weight}}
\font\gotic=eufm10
\newcommand\gots{{\hbox{\gotic\char'123}}}
\newenvironment{note}[1]{\par\addvspace{\medskipamount}\noindent
                         {\bf {#1}}\sl
                       }{\par\addvspace{\medskipamount}\rm}

\title{On Characters of Weyl Groups}

\author{Ron M. Adin%
\thanks{Department of Mathematics and Computer Science, Bar-Ilan University,
Ramat Gan 52900, Israel. Email: {\tt radin@math.biu.ac.il} }$^{\ \S}$
\and
Alexander Postnikov%
\thanks{Department of Applied Mathematics, Massachusetts
Institute of Technology, MA 02139, USA. Email: {\tt apost@math.mit.edu} }
\and Yuval Roichman%
\thanks{Department of Mathematics and Computer Science, Bar-Ilan University,
Ramat Gan 52900, Israel. Email: {\tt yuvalr@math.biu.ac.il} } 
\thanks{Research supported in part by  the Israel Science Foundation
and by internal research grants from Bar-Ilan University.}}
\date{(submitted: November 20, 1998;\ \ revised: April 30, 2000)}

\begin{document}

\maketitle

\begin {abstract}
In this note a combinatorial character formula related to
the symmetric group is generalized to an arbitrary finite Weyl group.
\end{abstract}

\section{The Case of the Symmetric Group}

The length $\ell(\pi)$ of a permutation
$\pi\in S_n$ is the number of inversions of $\pi$,
i.e., the number of pairs $(i,j)$ with $1\le i< j\le n$ and
$\pi(i)>\pi(j)$. 

For any permutation $\pi\in S_n$ let $m(\pi)$ be defined as
$$
m(\pi):=\cases
{(-1)^m,& if there exists $0\le m< n$ so that\cr
& $\pi(1)>\pi(2)>\dots>\pi(m+1)<\dots<\pi(n)$;\cr
0, & otherwise.\cr}
\leqno(1)
$$

Let $\mu=(\mu_1,\dots,\mu_t)$ be a partition of $n$, and let 
$S_\mu=S_{\mu_1}\times S_{\mu_2}\times\cdots\times S_{\mu_t}$
be the corresponding Young subgroup of $S_n$.
For any permutation
$\pi=r\cdot (\pi_1\times\cdots\times \pi_t)$, where
$\pi_i\in S_{\mu_i}$ $(1\le i\le t)$ and $r$
is a representative of minimal length of a left coset of $S_\mu$ in $S_n$,
define
$$
\weight_\mu(\pi):=\prod\limits_{i=1}^t m(\pi_i), \leqno(2)
$$
where $m(\pi_i)$ is defined in $S_{\mu_i}$ by (1).

\medskip

Denote by $\chi^k_\mu$ the value, at a conjugacy class of type $\mu$,
of the character of the natural $S_n$-action on the $k$-th homogeneous 
component of the coinvariant algebra.
The following combinatorial character formula was proved in 
[Ro2, Theorem 1].

\medskip

\begin{note}
\noindent{\bf Theorem.} With the above notations
$$
\chi^k_\mu=\suml_{\{\pi\in S_n:\ell(\pi)=k\}} \weight_\mu(\pi).
$$
\end{note}

\medskip

\section{Arbitrary Weyl Group}

Let $W$ be an arbitrary finite Weyl group.
Denote the set of positive roots by $\Phi_+$. 
Let $t_\alpha$ be the reflection corresponding to $\alpha\in \Phi_+$, and let 
$\check\alpha$ be the corresponding coroot.
Let $\alpha_i$
be the simple root corresponding to the simple reflection $s_i$.
 Denote by $\gots_w$ the Schubert
polynomial indexed by $w\in W$. 

The following theorem describes the action of the simple reflections
on the coinvariant algebra. This theorem is a 
reformulation of [BGG, Theorem 3.14 (iii)].

\medskip

\begin{note}
\noindent{\bf Theorem 1.}  For any simple reflection $s_i$ in $W$
and any $w\in W$,  
$$
s_i(\gots_w)=\cases
{\gots_w, & if $\ell(w s_i)>\ell(w)$;\cr
-\gots_w + \suml_{\{\alpha\in \Phi_+| \alpha\not=\alpha_i \wedge
 \ell(ws_it_{\alpha})=\ell(w)\}} \alpha_i(\check \alpha)
\gots_{ws_i t_\alpha}, &  if $\ell(w s_i)<\ell(w)$.\cr}
$$
\end{note}

\noindent{\bf Proof.} In the above notations, [BGG, Theorem 3.14 (iii)] states that
$$
s_i(\gots_w)=\cases
{\gots_w, & if $\ell(w s_i)>\ell(w)$;\cr
\gots_w - \suml_{\{\gamma\in \Phi_+| 
 \ell(t_\gamma ws_i)=\ell(w)\}} w (\alpha_i)(\check \gamma)
\gots_{t_\gamma ws_i}, &  if $\ell(w s_i)<\ell(w)$.\cr}
$$
Obviously, for any $\gamma\in \Phi_+$ there exists a unique
$\alpha\in \Phi_+$ such that $t_\gamma w s_i=ws_i t_\alpha$. In this case
$s_i w^{-1}(\gamma)=\alpha$. If $\ell(ws_i)<\ell(w)$ and $t_\gamma w s_i\not=w$
then the coefficient of $\gots_{ws_it_\alpha}=\gots_{t_\gamma w s_i}$  
in $s_i(\gots_w)$ is equal to
$$
-w(\alpha_i)(\check \gamma)= 
-\alpha_i({w^{-1} (\check \gamma)}) = 
-\alpha_i(s_i(\check\alpha)) =
\alpha_i(\check\alpha).
$$

If $t_\gamma w s_i=w$ then $t_\alpha=s_i$. 
Hence, the  coefficient of $\gots_w$ in $s_i(\gots_w)$ is
$1-w (\alpha_i)(\check \gamma)=1-\alpha_i(\check\alpha)=-1$ 
if $\ell(s_i w)<\ell(w)$, 
and $1$ otherwise.
\qed

\medskip

Let  $\langle\hbox{ },\hbox{ }\rangle$ be the inner product 
on the coinvariant algebra defined by
$\langle \gots_v,\gots_w\rangle=\delta_{v,w}$  (the Kronecker delta).
Theorem 1 implies

\medskip

\begin{note}
\noindent{\bf Corollary 2.} 
 Let $s_i$ be a simple reflection in  $W$,
 and let $z\in W$ such that
$\ell(z s_i)<\ell(z)$. Then for any $w\in W$ 
$$
\langle s_i (\gots_w),\gots_z\rangle=\cases
{0,& if $z\not=w$\cr
-1,& if $z=w$\cr} .
$$
\end{note}

\medskip

\noindent{\bf Proof.} For $z=w$ this follows from the second case of Theorem 1.
For $z\not= w$
if $\langle s_i (\gots_w),\gots_z\rangle\not=0$
then (by Theorem 1) $z=ws_it_\alpha$ for some $\alpha\in \Phi_+$ such that 
$\alpha\not=\alpha_i$ and
$\ell(ws_it_\alpha)=\ell(w)>\ell(ws_i)$. 
Now, for $\alpha\in \Phi_+$, 
$\ell(ws_it_\alpha)>\ell(ws_i)$ 
if and only if $ws_i(\alpha)\in \Phi_+$.
On the other hand, $\alpha_i\not= \alpha\in \Phi_+ \Rightarrow
s_i(\alpha)\in \Phi_+$. Since $ws_i(\alpha)\in \Phi_+$ it follows that
$\ell(wt_{s_i(\alpha)})>\ell(w)$. But $wt_{s_i(\alpha)}=ws_it_\alpha s_i$.
Hence, 
$\ell(zs_i)=\ell(ws_it_\alpha s_i)>\ell(w)=\ell(ws_it_\alpha)=\ell(z).$
\qed

\medskip

The following is, surprisingly, an exact Schubert analogue of 
a useful vanishing condition for Kazhdan-Lusztig coefficients
[Ro1, Lemma 4.3].

\medskip

\begin{note}
\noindent{\bf Corollary 3.} Let $s_i,s_j$ be
commuting  simple reflections in $W$, and
let $w,z\in W$ such that $\ell(ws_i)> \ell(w)$ and $\ell(zs_i)<\ell(z)$.
Then 
$$
\langle s_j(\gots_w), \gots_z \rangle= 0.
$$
\end{note}

\medskip

\noindent{\bf Proof.}  Obviously, $z\not=w$.
If $\ell(ws_j)>\ell(w)$ then our claim is an immediate consequence of Theorem 1.
Assume that $\ell(ws_j)<\ell(w)$, and
denote $\langle s_j (\gots_w), \gots_z \rangle$ by 
$b^{(j)}_z(w)$.
By Corollary 2
$$
s_i(1+s_j)(\gots_w)=
s_i \left(\suml_{\ell(zs_j)>\ell(z)} b^{(j)}_z(w) \gots_z\right)
=\suml_{\ell(zs_j)>\ell(z)} b^{(j)}_z(w) s_i (\gots_z).
$$
On the other hand, by Theorem 1, $\gots_w$
is an invariant under $s_i$. Thus,
$$
s_i(1+s_j)(\gots_w)=(1+s_j)s_i (\gots_w)= (1+s_j)(\gots_w)=
\suml_{\ell(zs_j)>\ell(z)} b^{(j)}_z(w) \gots_z.
$$
We conclude that
$$
\suml_{\ell(zs_j)>\ell(z)} b^{(j)}_z(w) (1-s_i)\gots_z=0.
$$
But
$$
\suml_{\ell(zs_j)>\ell(z)} b^{(j)}_z(w) (1-s_i)(\gots_z)=
\suml_{\ell(zs_j)>\ell(z)\wedge \ell(zs_i)<\ell(z)} b^{(j)}_z(w) (1-s_i)(\gots_z)=
$$
$$
=\suml_{\ell(zs_j)>\ell(z)\wedge  \ell(zs_i)<\ell(z)} b^{(j)}_z(w) [2\gots_z-\suml_{\ell(ts_i)>\ell(t)} b^{(i)}_t(z) \gots_t].
$$
This sum is equal to zero if and only if $b^{(j)}_z(w)=0$
for all $z$ with $\ell(zs_j)>\ell(z)$ and $\ell(zs_i)<\ell(z)$.

It remains to check the case in which $\ell(zs_j)<\ell(z)$.
By assumption $\ell(zs_i)<\ell(z)\Rightarrow z\not=w$.
 Corollary 2 completes the proof.
\qed

\medskip

 Let $H$ be a parabolic subgroup of $W$, which is 
isomorphic to a direct product of symmetric groups.
In the following definition we refer to cycle type and $\weight_\mu$
of elements in $H$ under the natural isomorhism, sending simple reflections of $H$ to simple reflections of $W$.

\medskip

\noindent{\bf Definition.} 
Let  $\mu$ be a cycle type of an element in $H$.
For any element 
$w=r\cdot \pi\in W$, where $\pi\in H$ and $r$ 
is the representative of minimal length of the left coset of $wH$ in $W$,
define
$$
\weight_\mu(w):=\weight_\mu(\pi). 
$$
Here $\weight_\mu(\pi)$ is defined as in Section 1.

\medskip

Note that  
$\weight_\mu$ is independent of the choice of $H$, provided that $H$
is isomorhic to a direct product of symmetric groups and that
 $\mu$ is the cycle type of some element in $H$.

Let $R^k$ be the $k$-th
homogeneous component of the coinvariant algebra of $W$. 
Denote by $\chi^k$ the $W$-character of 
$R^k$. Let $v_\mu\in H$ have cycle type $\mu$.
Then

\medskip

\begin{note}
\noindent{\bf Theorem 4.} With the above notations
$$
\chi^k(v_\mu)=\suml_{\{w\in W:\ell(w)=k\}} \weight_\mu(w).
$$
\end{note}

\noindent{\bf Proof.} Imitate the proof of [Ro2, Theorem 1].
Here Corollary 2 plays
 the role of [Ro2, Corollary 3.2] and implies an analogue of 
[Ro2, Corollary 3.3]. 
Alternatively, one can prove Theorem 4 by imitating the proof of 
 [Ro1, Theorems 1-2], where Corollary
3 plays the role of [Ro1, Lemma 4.3].
\qed

\medskip

\noindent{\bf Note:} A formally similar result appears also in Kazhdan-Lusztig theory.
The Kazhdan-Lusztig
characters of $W$ at $v_\mu$ may be represented as sums
of exactly the same weights, but over Kazhdan-Lusztig cells 
instead of Bruhat levels [Ro1, Corollary 3]. This curious analogy seems to deserve further
study.
For a $q$-analogue of the result for the symmetric group
see [APR]. A $q$-analogue of Theorem 4 is desired.

\end{document}